\documentclass[11pt,reqno]{amsart}
\usepackage{graphicx}
\usepackage{hyperref}
\usepackage{amssymb}
\usepackage{cite}
\usepackage{amsmath}
\usepackage{latexsym}
\usepackage{amscd}
\usepackage{amsthm}
\usepackage{mathrsfs}
\usepackage{url}

\newtheorem{thm}{Theorem}[section]

\newtheorem{lem}[thm]{Lemma}

\theoremstyle{definition}

\theoremstyle{remark}
\newtheorem{rem}[thm]{Remark}
\newtheorem*{ack}{Acknowledgment}
\def\R{\mathbb R}

\def\H{\mathbb H}
\def\SS{\mathbb S}

\def\f{\frac}

\def\ra{\rightarrow}
\def\pt{\partial}

\begin{document}
\title[Alexandrov-Fenchel type inequalities]{Alexandrov-Fenchel type inequalities for convex hypersurfaces in hyperbolic Space and in Sphere}
\author{Yong Wei}
\author{Changwei Xiong}
\address{Department of mathematical sciences, Tsinghua University, 100084, Beijing, P. R. China}
\email{wei-y09@mails.tsinghua.edu.cn}
\email{xiongcw10@mails.tsinghua.edu.cn}
\date{\today}

\maketitle

\begin{abstract}
In this paper, firstly, inspired by Nat\'{a}rio's recent work \cite{Na}, we use the isoperimetric inequality to derive some Alexandrov-Fenchel type inequalities for closed convex hypersurfaces in the hyperbolic space $\H^{n+1}$ and in the sphere $\SS^{n+1}$. We also get the rigidity in the spherical case. Secondly, we  use the inverse mean curvature flow in sphere \cite{gerh,Mak-Sch} to prove an optimal Sobolev type inequality for closed convex hypersurfaces in the sphere.
\end{abstract}

\section{Introduction}\label{sec:1}

Let $N^{n+1}(c)$ be the space form of constant sectional curvature $c$ and $\psi:\Sigma^n\ra N^{n+1}(c)$ be a closed hypersurface. Denote the $k$-th order mean curvature of $\Sigma$ by $p_k$ (see section \ref{sec:pre}). The inequalities about the integrals $\int_{\Sigma} p_k d\mu$ have been attracting many attentions for a long time. Among them the most famous one is the classical Minkowski inequality for closed convex surfaces $\Sigma\subset \R^3$, which can be written as
\begin{equation}\label{minkow-2}
\left(\f{1}{\omega_2}\int_{\Sigma} p_1 d\mu\right)^2\geq \f{|\Sigma|}{\omega_2},
\end{equation}
with equality if and only if $\Sigma$ is a sphere. Here $\omega_n$ is the area of $\SS^n(1)$ and $|\Sigma|=\int_{\Sigma}d\mu$ is the area of $\Sigma$ with respect to the induced metric from $\R^3$. The general inequality is the Alexandrov-Fenchel inequality \cite{Alex1,Alex2,Fe} which states that for convex hypersurfaces in the Euclidean space $\R^{n+1}$, there holds
\begin{eqnarray}
\frac 1{\omega_n}\int_{\Sigma}p_kd\mu &\geq& \left( \frac 1{\omega_n}\int_{\Sigma}p_ld\mu\right)^{\frac{n-k}{n-l}},\,0\leq l<k\leq n,
\end{eqnarray}
with equality if and only if $\Sigma$ is a sphere. See \cite{ben,CW,GL,McC,Schn} for other references on  Alexandrov-Fenchel inequalities for closed hypersurfaces in Euclidean space $\R^{n+1}$,

It is natural to generalize the Minkowski inequality and Alexandrov-Fenchel inequalities to the hypersurfaces in the hyperbolic space $\H^{n+1}$. See for example \cite{BM,GS,Na}. Recently, the optimal Alexandrov-Fenchel type inequalities in $\H^{n+1}$ (see \cite{GWW1,GWW2,LWX,WX}) were obtained as following: For $1\leq k\leq n$ and closed horospherical convex hypersurface $\Sigma\subset\H^{n+1}$, there holds
\begin{eqnarray}\label{alex-hyper}
\frac 1{\omega_n}\int_{\Sigma}p_kd\mu &\geq& \left( \left(\f{|\Sigma|}{\omega_n}\right)^{\frac 2k}+\left(\f{|\Sigma|}{\omega_n}\right)^{\frac{2(n-k)}{kn}}\right)^{\frac k2},
\end{eqnarray}
with equality if and only if $\Sigma$ is a geodesic sphere in $\H^{n+1}$. In particular, when $k=2$, Li-Wei-Xiong \cite{LWX} proved that \eqref{alex-hyper} holds under a weaker condition that $\Sigma$ is star-shaped and 2-convex. In the proof of \eqref{alex-hyper}, the geometric flow was used and was an important tool. There are other papers on weighted Alexandrov-Fenchel type inequalities for hypersurfaces in the hyperbolic space $\H^{n+1}$, which were also proved by using the geometric flow. See \cite{BHW,de-girao,GWW3}.

\vskip 2mm
The Minkowski inequality and the Alexandrov-Fenchel inequalities can be viewed as the generalizations of the classical isoperimetric inequality, which compares the area of the hypersurface $\Sigma$ and the volume of the enclosed domain by $\Sigma$. The Minkowski inequality \eqref{minkow-2} was used by Minkowski himself to prove the isoperimetric inequality for closed convex surfaces (see \cite{Mink,Oss} ). Recently, J.Nat\'{a}rio \cite{Na} reversed Minkowski's idea, and derived a new Minkowski-type inequality for closed convex surfaces in the hyperbolic space $\H^3$ by using the isoperimetric inequality. In this paper, firstly, we deal with the higher dimensional case by adapting Nat\'{a}rio's method \cite{Na}. We will derive some new Alexandrov-Fenchel type inequalities for closed convex hypersurfaces in $\H^{n+1}$ and in $\SS^{n+1}$, starting from the isoperimetric inequality.

Let $\Sigma$ be a closed convex hypersurface in $\H^{n+1}$. Then $\Sigma$ is a smooth boundary of a convex body  $\Omega$. Recall that a subset $\Omega$ of a Riemannian manifold is called a convex body if it is convex with nonempty interior. The subset $\Omega$ is convex means that for any two points $x_0,x_1\in\Omega$, there exists a minimizing geodesic connecting $x_0$ and $x_1$ which are contained in $\Omega$. Inspired by Nat\'{a}rio's paper \cite{Na}, we define a family of parallel hypersurfaces $\Sigma_t=\psi_t(\Sigma)$, where $\psi_t(x)=\exp_{\psi(x)}(t\nu(x)), x\in\Sigma$,$\nu(x)$ is the outward unit normal of $\Sigma$. Since the expression of the geodesic in space forms can be written down explicitly (see \cite{MR}), we can directly compute the area of $\Sigma_t$ and the volume of the domain $\Omega_t$ enclosed by $\Sigma_t$. Note that in Nat\'{a}rio's paper, the area of $\Sigma_t$ was obtained by using the first and second variation formula with the help of the Gauss-Bonnet formula. Since $\H^{n+1}$ has constant negative curvature $c=-1$ and $\Sigma$ is convex, it follows that the parallel hypersurfaces $\Sigma_t$ can be well defined for all $t\geq 0$ (see \cite{Na}). Define a function $r(t)$ such that $|\Sigma_t|=|S_{r(t)}|$, then the isoperimetric inequality (see \cite{S,Ros}) implies that $Vol(\Omega_t)\leq Vol(B_{r(t)})$, where $S_{r(t)}$ and $B_{r(t)}$ are the geodesic sphere and geodesic ball of radius $r(t)$, respectively. Then applying the isoperimetric inequality to $\Sigma_t$ for sufficiently large $t$, we obtain the following Alexandrov-Fenchel type inequalities in $\H^{n+1}$.

\begin{thm}\label{thm1}
Let $\Sigma^n$ be a closed convex  hypersurface in $\H^{n+1}$ with $n\geq 3$. Then there holds
\begin{eqnarray}\label{eq-thm1}
\sum_{k=0}^n \f{2k-n}{n\omega_n} \int_{\Sigma} C_n^k p_k d\mu &\geq&  \left(\f{1}{\omega_n}\sum_{k=0}^n\int_{\Sigma} C_n^k p_k d\mu\right)^{\f{n-2}{n}}.
\end{eqnarray}
\end{thm}

A direct calculation shows that if $\Sigma$ is a geodesic sphere, then the equality in \eqref{eq-thm1} holds. However, we do not obtain the rigidity (i.e.,we don't know whether the equality in \eqref{eq-thm1} implies that $\Sigma$ is a geodesic sphere). At the end of section \ref{sec:3}, we will give a remark that when the hypersurface $\Sigma\subset\H^{n+1}$ is sufficiently small, the inequality \eqref{eq-thm1} reduces to one of the Alexandrov-Fenchel inequalities in Euclidean space.

\vskip 2mm
Next we will use the same method to derive inequalities for closed convex  hypersurfaces in $\SS^{n+1}$. In this case, we can prove the rigidity result.

\begin{thm}\label{thm2}
Let $\Sigma^n$ be a closed convex hypersurface in $\SS^{n+1}$ with $n\geq 2$. Then there holds
\begin{eqnarray}\label{eq-thm2}
\omega_n &\leq& \sum_{s=\f{1-(-1)^n}{2},+2}^n \sqrt{(E(s))^2+(F(s))^2},
\end{eqnarray}
where
\begin{align*}
E(s)&=\sum_{p+q=\f{n+s}{2}or\f{n-s}{2}}\sum_{q\leq k\leq n-p,2|k}C_n^k \f{1}{2^n}C_{n-k}^p C_k^q (-1)^{\f{k}{2}+k-q}\int_{\Sigma} p_k d\mu,\\
F(s)&=\sum_{p+q=\f{n+s}{2}or\f{n-s}{2}}\sum_{q\leq k\leq n-p,2\nmid k}C_n^k \f{1}{2^n}C_{n-k}^p C_k^q (-1)^{[\f{k}{2}]+k-q}(-1)^{\chi_{\{2(p+q)-n\leq 0\}}}\int_{\Sigma} p_k d\mu,
\end{align*}
and ``$+2$'' means that the step length of the summation for $s$ is $2$.
Moreover, the equality holds in \eqref{eq-thm2} if and only if $\Sigma^n$ is a geodesic sphere.
\end{thm}

When $n=2$, it is easy to check that
\begin{align*}
E(0)&=2\pi, & F(0)&=0,\\
E(2)&=|\Sigma|-2\pi, & F(2)&=\int_{\Sigma} p_1 d\mu,
\end{align*}
where the Gauss-Bonnet Theorem $|\Sigma|+\int_{\Sigma} p_2 d\mu=4\pi$ (see section \ref{sec:L_k}) was used. So \eqref{eq-thm2} implies the Minkowski-type inequality in the sphere
\begin{eqnarray}\label{eq-thm2-1}
\left(\int_{\Sigma} p_1 d\mu\right)^2&\geq& |\Sigma|(4\pi-|\Sigma|),
\end{eqnarray}
which is just the Theorem 0.2 in \cite{Na}. See also \cite{Bla,Kno,Sant}. Note that Makowski and Scheuer \cite{Mak-Sch} recently also proved \eqref{eq-thm2-1} by using a different method involving the inverse curvature flow in sphere. To get more feeling of the inequality \eqref{eq-thm2}, we also give the precise expressions of \eqref{eq-thm2} in the case of $n=3$ and $n=4$. See Remark \ref{rem:5.1} in Section \ref{sec:4}.

\vskip 2mm
Finally, in the last part of this paper, we use the inverse mean curvature flow in the sphere \cite{Mak-Sch,gerh} to prove the following optimal inequalities for convex hypersurfaces in sphere $\SS^{n+1}$.
\begin{thm}\label{alex-fenchel-sphere}
Let $\Sigma^n$ be a closed and strictly convex hypersurface in $\SS^{n+1}$. Then we have the following optimal inequalities ($k\leq n/2$)
\begin{eqnarray}\label{p_k-sphere}
  \int_{\Sigma}L_kd\mu&\geq& C_{n}^{2k}(2k)!\omega_{n}^{\frac {2k}{n}}|\Sigma|^{\frac{n-2k}{n}}.
\end{eqnarray}
Equality holds in \eqref{p_k-sphere} if and only if $\Sigma$ is a geodesic sphere. Here $L_k$ is the Gauss-Bonnet curvature of the induced metric on $\Sigma$ (see section \ref{sec:L_k} for detail).
\end{thm}

The proof of Theorem \ref{alex-fenchel-sphere} is similar with the proof in \cite{BHW,de-girao,GL,GWW1,GWW2,LWX}. We define a curvature quantity $Q(t)$ which is monotone non-increasing under the inverse mean curvature flow in the sphere. Then we obtain the inequality \eqref{p_k-sphere} by comparing the initial value $Q(0)$ with the limit $\lim_{t\ra T^*}Q(t)$. We remark that since $\Sigma$ is a closed convex hypersurface in $\SS^{n+1}$, a well-known result due to do Carmo and Warner \cite{do-carmo} implies that $\Sigma$ is embedded, homeomorphic to $n$-sphere and is contained in an open hemisphere.

When $k=1$, Theorem \ref{alex-fenchel-sphere} says
\begin{eqnarray}\label{p_2-sphere}
\int_{\Sigma}p_2d\mu+|\Sigma| &\geq& \omega_n^{\frac 2n}|\Sigma|^{\frac{n-2}n} ,
\end{eqnarray}
which was already proved by Makowski and Scheuer in \cite{Mak-Sch}. One can compare \eqref{p_2-sphere} with the case $k=2$ of the Alexandrov-Fenchel type inequality \eqref{alex-hyper} in $\H^{n+1}$, i.e.,
\begin{eqnarray}\label{p_2-hyperbolic}
\int_{\Sigma}p_2d\mu-|\Sigma| &\geq& \omega_n^{\frac 2n}|\Sigma|^{\frac{n-2}n} ,
\end{eqnarray}
which was proved by Li-Wei-Xiong \cite{LWX} for star-shaped and 2-convex hypersurfaces in $\H^{n+1}$. Note that both \eqref{p_2-sphere} and \eqref{p_2-hyperbolic} have the rigidity results that the equality holds if and only if $\Sigma$ is a geodesic hypersphere.  For $k>1$, the inequalities of the same type as \eqref{p_k-sphere} were proved by Ge-Wang-Wu \cite{GWW1,GWW2} for horospherical convex hypersurfaces in the hyperbolic space $\H^{n+1}$.

\begin{ack}
The authors would like to thank Professor Haizhong Li for constant encouragement and help. The research of the authors was supported by NSFC No. 11271214.
\end{ack}

\section{Prelimilaries}\label{sec:pre}

\subsection{$k$-th order mean curvature}
Let $\Sigma$ be a closed hypersurface in $N^{n+1}(c)$ with unit outward normal $\nu$. The second fundamental form $h$ of $\Sigma$ is defined by
\begin{equation*}
    h(X,Y)=\langle\bar{\nabla}_X\nu,Y\rangle
\end{equation*}
for any two tangent fields $X,Y$. For an orthonormal basis $\{e_1,\cdots,e_n\}$ of $\Sigma$, the components of the second fundamental form are given by $h_{ij}=h(e_i,e_j)$ and $h_i^j=g^{jk}h_{ki}$, where $g$ is the induced metric on $\Sigma$. The principal curvature $\kappa=(\kappa_1,\cdots,\kappa_{n})$ are the eigenvalues of $h$ with respect to $g$. The $k$-th order mean curvature of $\Sigma$ for $1\leq k\leq n$ are defined as
\begin{equation}
  p_k=\frac 1{C_n^k}\sigma_k(\kappa)=\frac 1{C_n^k}\sum_{i_1<i_2<\cdots<i_k}\kappa_{i_1}\cdots\kappa_{i_k},
\end{equation}
or equivalently as
\begin{equation}
  p_k=\frac 1{C_n^k}\sigma_k(h_i^j)=\frac 1{C_n^kk!}\delta_{j_1\cdots j_k}^{i_1\cdots i_k}h_{i_1}^{j_1}\cdots h_{i_k}^{j_k},
\end{equation}
where $\delta_{j_1\cdots j_k}^{i_1\cdots i_k}$ is the generalized Kronecker delta given by
\begin{eqnarray*}
  \delta_{j_1\cdots j_k}^{i_1\cdots i_k} &=& \det\left(
                                                   \begin{array}{cccc}
                                                     \delta_{j_1}^{i_1} & \delta_{j_1}^{i_2} & \cdots & \delta_{j_1}^{i_k} \\
                                                     \delta_{j_2}^{i_1} & \delta_{j_2}^{i_2} & \cdots & \delta_{j_2}^{i_k} \\
                                                     \vdots & \vdots & \vdots & \vdots \\
                                                     \delta_{j_k}^{i_1} & \cdots & \cdots & \delta_{j_k}^{i_k} \\
                                                   \end{array}
                                                 \right)
\end{eqnarray*}
We have the following Newton-MacLaurin inequalities (see, e.g.,\cite{guan}).
\begin{lem}\label{newton_ineq}
For $\kappa\in \overline{\Gamma}^+_k$, where $\overline{\Gamma}^+_k$ is the closure of the Garding cone
\begin{equation*}
  \Gamma_k^+=\{\kappa\in\R^n|~p_j(\kappa)>0,~\forall j\leq k\},
\end{equation*}
we have the following Newton-MacLaurin inequalities
\begin{align*}
    &p_1p_{k-1}\geq p_k\\
    &p_1\geq p_2^{1/2}\geq\cdots\geq p_k^{1/k}.
\end{align*}
Equalities hold if and only if $\kappa=c(1,\cdots,1)$ for some constant $c$.
\end{lem}

\subsection{Gauss-Bonnet curvature $L_k$} \label{sec:L_k}
Given an n-dimensional Riemannian manifold $(M,g)$, the Gauss-Bonnet curvature $L_k (k\leq n/2)$ is defined by (see, e.g., \cite{GWW0,GWW2})
\begin{equation}\label{def-l_k}
  L_k=\frac 1{2^k}\delta_{j_1j_2\cdots j_{2k-1}j_{2k}}^{i_1i_2\cdots i_{2k-1}i_{2k}}R_{i_1i_2}^{\quad j_1j_2}\cdots R_{i_{2k-1}i_{2k}}^{\qquad j_{2k-1}j_{2k}}.
\end{equation}

For a closed hypersurface $\Sigma^{n}\subset\R^{n+1}$, recall the Gauss equation
\begin{equation*}
  R_{ij}^{~~~kl}=h_i^kh_j^l-h_i^lh_j^k.
\end{equation*}
Then the Gauss-Bonnet curvature of the induced metric on $\Sigma^{n}\subset\R^{n+1}$ is
\begin{eqnarray}
  L_k&=&\delta_{j_1j_2\cdots j_{2k-1}j_{2k}}^{i_1i_2\cdots i_{2k-1}i_{2k}}h_{i_1}^{j_1}h_{i_2}^{j_2}\cdots h_{i_{2k-1}}^{j_{2k-1}}h_{i_{2k}}^{j_{2k}}\nonumber\\
  &=&(2k)!C_n^{2k}p_{2k}.\label{L_k-R}
\end{eqnarray}
For a closed hypersurface $\Sigma^{n}\subset\SS^{n+1}$, the Gauss equations are
\begin{equation*}
  R_{ij}^{~~~kl}=(h_i^kh_j^l-h_i^lh_j^k)+(\delta_i^k\delta_j^l-\delta_i^l\delta_j^k).
\end{equation*}
Then by a straightforward calculation, we have
\begin{eqnarray}
L_k &=& \delta_{j_1j_2\cdots j_{2k-1}j_{2k}}^{i_1i_2\cdots i_{2k-1}i_{2k}}(h_{i_1}^{j_1}h_{i_2}^{j_2}+\delta_{i_1}^{j_2}\delta_{i_2}^{j_2})\cdots (h_{i_{2k-1}}^{j_{2k-1}}h_{i_{2k}}^{j_{2k}}+\delta_{i_{2k-1}}^{j_{2k-1}}\delta_{i_{2k}}^{j_{2k}}) \nonumber\\
   &=& \sum_{i=0}^kC_k^i(n-2k+1)(n-2k+2)\cdots (n-2k+2i)(2k-2i)!\sigma_{2k-2i} \nonumber\\
  &=& \sum_{i=0}^kC_k^i(n-2k+1)(n-2k+2)\cdots (n-2k+2i)(2k-2i)!\nonumber\\
  &&\quad \times C_{n}^{2k-2i}p_{2k-2i}\nonumber\\
  &=&\sum_{i=0}^kC_k^i\frac {n!}{(n-2k)!}p_{2k-2i}\nonumber\\
  &=&C_{n}^{2k}(2k)!\sum_{i=0}^kC_k^ip_{2k-2i}.
\end{eqnarray}
Similarly, for a closed hypersurface $\Sigma^{n}\subset\H^{n+1}$, its Gauss-Bonnet curvature is
\begin{equation}
 L_k=C_{n}^{2k}(2k)!\sum_{i=0}^kC_k^i(-1)^ip_{2k-2i}.
\end{equation}

Finally, note that through our paper, we assume that the hypersurface $\Sigma\subset N^{n+1}(c)$ is closed and convex. It follows that $\Sigma$ is homeomorphic to the n-sphere.  Then if the dimension of $\Sigma$ is even, the Gauss-Bonnet-Chern theorem \cite{Chern} implies that
\begin{equation}\label{gbc-thm}
  \int_{\Sigma}L_{\frac n2}d\mu=n!\omega_n.
\end{equation}
\eqref{gbc-thm} will be used in the following sections. Also \eqref{gbc-thm} shows that when $2k=n$, the inequality \eqref{p_k-sphere} is an equality.

\section{The Euclidean case}\label{sec:2}

To demonstrate the method which will be used to prove Theorem \ref{thm1} and Theorem \ref{thm2}, in this section, we first consider the simple case that $N^{n+1}(c)=\R^{n+1}$.  As said in section \ref{sec:1}, we define a family of parallel hypersurfaces $\Sigma_t=\psi_t(\Sigma)$, where $\psi_t(x)=\exp_{\psi(x)}(t\nu(x)), x\in\Sigma$. In this case, $\psi_t=\psi+t\nu$ (see \cite{MR}) and so $(\psi_t)_*(e_i)=(1+t\kappa_i)e_i$. Therefore the area element of $\Sigma_t$ is
\begin{equation*}
d\mu_t=(1+t\kappa_1)\cdots(1+t\kappa_n)d\mu,
\end{equation*}
which implies that the area $|\Sigma_t|$ of the parallel hypersurfaces $\Sigma_t$ are equal to
\begin{align*}
|\Sigma_t|&=\int_{\Sigma} (1+t\kappa_1)\cdots(1+t\kappa_n)d\mu\\
     &=\sum_{k=0}^n\int_{\Sigma}C_n^kp_kd\mu \,t^k.
\end{align*}
For abbreviation, we denote by $|V|$ the volume of the domain $\Omega$ enclosed by $\Sigma$, and $|V_t|$ the volume of the domain $\Omega_t$ enclosed by $\Sigma_t$. By integrating and noting that $\Sigma_t=\psi_t(\Sigma)$ are parallel hypersurfaces of $\Sigma$ given by $\psi_t(x)=\exp_{\psi(x)}(t\nu(x))$ for $x\in\Sigma$, we can obtain
\begin{align*}
|V_t|&=|V|+\int_0^t|\Sigma_s|ds\\
     &=|V|+\sum_{k=0}^n\int_{\Sigma}C_n^kp_kd\mu\frac 1{k+1}t^{k+1}.
\end{align*}

Now since the initial hypersurface is convex, the $\psi_t$ is well defined for all $t>0$. For $\forall t\geq 0$, the isoperimetric inequality (see \cite{Oss}) in Euclidean space $\R^{n+1}$ implies
\begin{equation}\label{isop-eucl}
\left(\f{|\Sigma_t|}{\omega_n}\right)^{n+1}\geq ((n+1)\f{|V_t|}{\omega_n})^n.
\end{equation}

If $n$ is odd, then comparing the coefficient of $t^{n(n+1)}$ in \eqref{isop-eucl} yields
\begin{equation}\label{eq-euclid-odd}
\int_{\Sigma} p_nd\mu\geq \omega_n,
\end{equation}
which is a special Alexandrov-Fenchel inequality.

If $n$ is even, \eqref{L_k-R} and the Gauss-Bonnet-Chern Theorem \eqref{gbc-thm} imply that $ \int_{\Sigma} p_n d\mu =\omega_n$. Thus expanding the two sides of \eqref{isop-eucl} and comparing the coefficients of $t^{n(n+1)}$, $t^{n(n+1)-1}$ and $t^{n(n+1)-2}$, we can get
\begin{eqnarray}\label{eq-euclid}
\left(\f{1}{\omega_n}\int_{\Sigma} p_{n-1}d\mu\right)^2&\geq& \f{1}{\omega_n} \int_{\Sigma} p_{n-2}d\mu,
\end{eqnarray}
which is also an Alexandrov-Fenchel inequality. In particular, when $n=2$, \eqref{eq-euclid} reduces to the classical Minkowski inequality \eqref{minkow-2}.

\section{The Hyperbolic case}\label{sec:3}

Now assume $N^{n+1}(c)=\H^{n+1}$. Since the initial hypersurface is convex,  the parallel hypersurfaces $\Sigma_t$ is well defined for all $t\geq 0$ (\cite{Na}). Note that $\psi_t=\cosh t \psi +\sinh t \nu$ (see \cite{MR}) and so $(\psi_t)_*(e_i)=(\cosh t +\sinh t \kappa_i)e_i$. Therefore the area element of $\Sigma_t$ is
\begin{equation*}
d\mu_t=(\cosh t +\sinh t \kappa_1)\cdots(\cosh t +\sinh t \kappa_n)d\mu,
\end{equation*}
which implies
\begin{align*}
|\Sigma_t|&=\int_{\Sigma} (\cosh t +\sinh t \kappa_1)\cdots(\cosh t +\sinh t \kappa_n)d\mu\\
     &=\sum_{k=0}^n\int_{\Sigma}C_n^kp_kd\mu\cosh^{n-k}t\sinh^kt.
\end{align*}
Then by integrating we can obtain
\begin{align*}
|V_t|&=|V|+\int_0^t|\Sigma_s|ds\\
     &=|V|+\sum_{k=0}^n\int_{\Sigma}C_n^kp_kd\mu\int_0^t\cosh^{n-k}s\sinh^ksds.
\end{align*}

On the other hand, it is well known that the area of a geodesic sphere $S_r$ and the volume of a geodesic ball $B_r$ with radius $r$ in the hyperbolic space $\H^{n+1}$ are
\begin{align*}
&S(r):=|S_r|=\omega_n\sinh^n r ,\\
&V(r):=Vol(B_r)=\omega_n\int_0^r \sinh^n s ds.
\end{align*}
Now define a function $r(t)$ such that $|\Sigma_t|=S(r(t))$. That is,
\begin{equation}\label{eq0}
\sum_{k=0}^n \cosh^{n-k} t \sinh^k t \cdot \int_{\Sigma} C_n^k p_k d\mu=\omega_n \sinh^n r(t).
\end{equation}
Then the isoperimetric inequality (see \cite{S,Ros}) implies
\begin{equation}\label{eq1}
|V_t|\leq V(r(t)), \forall t \geq 0.
\end{equation}
From this inequality, we can get some information for $\Sigma$.

First we will get a rough estimate for $r(t)$. When $t\rightarrow +\infty$, $\cosh^{n-k} t\sinh^k t=\sinh^n t(1+o(1))$. Thus from $|\Sigma_t|=S(r(t))$, we get
\begin{equation*}
\sinh^n t(1+o(1))\sum_{k=0}^n\int_{\Sigma} C_n^k p_k d\mu=\omega_n \sinh^n r(t),
\end{equation*}
which implies
\begin{equation}\label{estimate}
r(t)=t+\f{1}{n}\ln (\f{1}{\omega_n}\sum_{k=0}^n\int_{\Sigma} C_n^k p_k d\mu)+o(1).
\end{equation}

However, this estimate for $r(t)$ is not enough. For our purpose we should make better use of $|\Sigma_t|=S(r(t))$ as follows. The case of $n=2$ was considered by Nat\'{a}rio in \cite{Na}, so we assume that $n\geq 3$ in the following calculation. Since we will examine \eqref{eq1} for sufficiently large $t$, we only care about the terms involving $e^{nt}$ and $e^{(n-2)t}$. The other terms are $o(e^{(n-2)t})$. It is straightforward to check that
\begin{equation*}
\cosh^{n-k} t \sinh^k t=\f{1}{2^n}e^{nt}+\f{1}{2^n}(n-2k)e^{(n-2)t}+\cdots.
\end{equation*}
Consequently \eqref{eq0} implies
\begin{align}
&\f{1}{2^n}\sum_{k=0}^n(e^{nt}+(n-2k)e^{(n-2)t}+\cdots)\int_{\Sigma} C_n^k p_k d\mu \nonumber\\
&=\omega_n(\f{1}{2^n}e^{nr}-\f{1}{2^n}n e^{(n-2)r}+\cdots).\label{eq11}
\end{align}

On the other hand,
\begin{align*}
|V_t|&=|V|+\sum_{k=0}^n \int_0^t \cosh^{n-k} s \sinh^k s ds\cdot \int_{\Sigma} C_n^k p_k d\mu\\
     &=|V|+\f{1}{2^n}\sum_{k=0}^n\int_0^t (e^{ns}+(n-2k)e^{(n-2)s}+\cdots) ds \int_{\Sigma} C_n^k p_k d\mu\\
     &=|V|+\f{1}{2^n}\sum_{k=0}^n(\f{1}{n} e^{nt}+\f{n-2k}{n-2}e^{(n-2)t}+\cdots) \int_{\Sigma} C_n^k p_k d\mu\\
     &=\f{1}{2^n}\f{1}{n}e^{nt} \sum_{k=0}^n \int_{\Sigma} C_n^k p_k d\mu+\f{1}{2^n} e^{(n-2)t}\sum_{k=0}^n \f{n-2k}{n-2}\int_{\Sigma} C_n^k p_k d\mu+\cdots,
\end{align*}
and
\begin{align*}
V(r(t))&=\omega_n \int_0^r \sinh^n s ds\\
       &=\omega_n \int_0^r (\f{1}{2^n}e^{ns}-\f{1}{2^n}n e^{(n-2)s}+\cdots) ds\\
       &=\f{\omega_n}{2^n}\f{1}{n} e^{nr}-\f{\omega_n}{2^n}\f{n}{n-2} e^{(n-2)r}+\cdots.
\end{align*}
Now taking \eqref{eq11} into account yields
\begin{align*}
V(r(t))&=\f{\omega_n}{2^n} e^{(n-2)r}+\f{1}{2^n}\f{1}{n}\sum_{k=0}^n(e^{nt}+(n-2k)e^{(n-2)t}+\cdots)\int_{\Sigma} C_n^k p_k d\mu\\
       &-\f{\omega_n}{2^n}\f{n}{n-2} e^{(n-2)r}\\
       &=\f{\omega_n}{2^n}(-\f{2}{n-2})e^{(n-2)r}+\f{1}{2^n}\f{1}{n}e^{nt}\sum_{k=0}^n \int_{\Sigma} C_n^k p_k d\mu\\
       &+\f{1}{2^n}\f{1}{n}e^{(n-2)t}\sum_{k=0}^n (n-2k)\int_{\Sigma} C_n^k p_k d\mu+\cdots.
\end{align*}

Noting \eqref{estimate}, we have
\begin{align*}
V(r(t))&=\f{\omega_n}{2^n}(-\f{2}{n-2})e^{(n-2)t}(\f{1}{\omega_n}\sum_{k=0}^n\int_{\Sigma} C_n^k p_k d\mu)^{\f{n-2}{n}}\\
       &+\f{1}{2^n}\f{1}{n}e^{nt}\sum_{k=0}^n \int_{\Sigma} C_n^k p_k d\mu+\f{1}{2^n}\f{1}{n}e^{(n-2)t}\sum_{k=0}^n (n-2k)\int_{\Sigma} C_n^k p_k d\mu+\cdots.
\end{align*}

Now $|V_t|\leq V(r(t)), t\rightarrow +\infty$ gives us
\begin{equation*}
\f{1}{2^n}\sum_{k=0}^n (\f{n-2k}{n-2}-\f{n-2k}{n})\int_{\Sigma} C_n^k p_k d\mu \leq \f{\omega_n}{2^n}(-\f{2}{n-2})\left(\f{1}{\omega_n}\sum_{k=0}^n\int_{\Sigma} C_n^k p_k d\mu\right)^{\f{n-2}{n}},
\end{equation*}
or equivalently
\begin{equation}\label{eq-hyper}
\sum_{k=0}^n \f{2k-n}{n} \int_{\Sigma} C_n^k p_k d\mu \geq \omega_n \left(\f{1}{\omega_n}\sum_{k=0}^n\int_{\Sigma} C_n^k p_k d\mu\right)^{\f{n-2}{n}}, n\geq 3.
\end{equation}
Hence we complete the proof of Theorem \ref{thm1}.
\begin{rem}
It is easy to check that for a geodesic sphere in $\H^{n+1}$, the equality holds in \eqref{eq-hyper}. However this method can not yield the rigidity result, i.e., we can not conclude that $\Sigma$ is a geodesic sphere if the equality holds in \eqref{eq-hyper}.
\end{rem}
\begin{rem}\label{rem:hyper}
We also remark that for small hypersurface $\Sigma\subset\H^{n+1}$ (i.e., with small diameter), the inequality \eqref{eq-hyper} can reduce to the Euclidean inequalities \eqref{eq-euclid-odd} and \eqref{eq-euclid}. For example, we first assume $n=4$. For 4-dimensional hypersurface $\Sigma\subset\H^5$, the Gauss-Bonnet-Chern formula \eqref{gbc-thm} implies
\begin{equation}\label{gb-hyper}
  \int_{\Sigma}(p_4-2p_2+1)d\mu=\frac 1{4!}\int_{\Sigma}L_2d\mu=\omega_4.
\end{equation}
Substituting \eqref{gb-hyper} into the inequality \eqref{eq-hyper} gives that
\begin{equation*}
  \left(1+\frac 2{\omega_4}\int_{\Sigma}(p_3+p_2-p_1-1)d\mu\right)^2\geq 1+\frac 4{\omega_4}\int_{\Sigma}(p_3+2p_2+p_1)d\mu.
\end{equation*}
Expanding the left hand side of the above inequality, and comparing both sides by orders (note that $\Sigma$ is a small hypersurface), we obtain that
\begin{eqnarray}
\left(\frac 1{\omega_4}\int_{\Sigma}p_3d\mu\right)^2 &\geq& \frac 1{\omega_4}\int_{\Sigma}p_2d\mu.
\end{eqnarray}
This is just the inequality \eqref{eq-euclid} for hypersurface in Euclidean space $\R^5$. For general even dimensional case, by using the Gauss-Bonnet-Chern formula,
\begin{equation*}
  \int_{\Sigma}\sum_{k=0}^{n/2}C_{n/2}^k(-1)^kp_{n-2k}d\mu=\frac 1{n!}\int_{\Sigma}L_{\frac n2}d\mu=\omega_n,
\end{equation*}
we can also reduce the inequality \eqref{eq-hyper} to the Euclidean version \eqref{eq-euclid} for small hypersurface $\Sigma\subset\H^{n+1}$. For odd dimensional case, the argument is similar.
\end{rem}

\section{The Spherical case}\label{sec:4}

Now assume $N^{n+1}(c)=\SS^{n+1}$. Since the initial hypersurface is convex, $\psi_t$ is well defined for $t\in [0,\f{\pi}{2})$ (\cite{Na}). Note that $\psi_t=\cos t \psi +\sin t \nu$($t\in [0,\f{\pi}{2})$) (see \cite{MR}) and so $(\psi_t)_*(e_i)=(\cos t +\sin t \kappa_i)e_i$. Therefore the area element of $\Sigma_t$ is
\begin{equation*}
d\mu_t=(\cos t +\sin t \kappa_1)\cdots(\cos t +\sin t \kappa_n)d\mu,
\end{equation*}
which implies
\begin{align*}
|\Sigma_t|&=\int_{\Sigma} (\cos t+\sin t \kappa_1)\cdots(\cos t+\sin t \kappa_n)d\mu\\
     &=\sum_{k=0}^n\int_{\Sigma}C_n^kp_kd\mu\cos^{n-k}t\sin^kt.
\end{align*}
Then by integrating we can also obtain
\begin{align*}
|V_t|&=|V|+\int_0^t|\Sigma_s|ds\\
     &=|V|+\sum_{k=0}^n\int_{\Sigma}C_n^kp_kd\mu\int_0^t\cos^{n-k}s\sin^ksds.
\end{align*}

On the other hand, it is well known that the area of a geodesic sphere $S_r$ and the volume of a geodesic ball $B_r$ with radius $r$ in the sphere $\SS^{n+1}$ are
\begin{align*}
S(r)&=\omega_n\sin^n r ,\\
V(r)&=\omega_n\int_0^r \sin^n s ds.
\end{align*}
Now since $|V_t|$ is increasing in $t$, when $t$ satisfies $|V_t|=V(\f{\pi}{2})=\omega_n\int_0^\f{\pi}{2} \sin^n r dr $, the isoperimetric inequality (see \cite{Ros}) implies $|\Sigma_t|\geq S(\f{\pi}{2})=\omega_n$ for this $t$. Therefore, a weaker requirement is
\begin{equation}\label{eq:Sphere}
\max_{t\in [0,\f{\pi}{2})} |\Sigma_t|\geq \omega_n.
\end{equation}

Then the key point is to estimate the $\max_{t\in [0,\f{\pi}{2})} |\Sigma_t|$. Direct computation shows that
\begin{align*}
\cos^{n-k} t \sin^k t&=(\f{e^{it}+e^{-it}}{2})^{n-k}(\f{e^{it}-e^{-it}}{2i})^{k}\\
                     &=\f{1}{2^n}\sum_{p=0}^{n-k}\sum_{q=0}^{k}C_{n-k}^p C_k^q \cos((2(p+q)-n)t-\f{k}{2}\pi)(-1)^{k-q}.
\end{align*}
Then
\begin{align*}
|\Sigma_t|&=\sum_{k=0}^n C_n^k \cos^{n-k} t \sin^k t \int_{\Sigma} p_k d\mu\\
     &=\sum_{k=0}^n C_n^k \f{1}{2^n}\sum_{p=0}^{n-k}\sum_{q=0}^{k}C_{n-k}^p C_k^q \cos((2(p+q)-n)t-\f{k}{2}\pi)(-1)^{k-q}\int_{\Sigma} p_k d\mu\\
     &=\sum_{0\leq k\leq n,2|k}C_n^k \f{1}{2^n}\sum_{p=0}^{n-k}\sum_{q=0}^{k}C_{n-k}^p C_k^q (-1)^{\f{k}{2}}\cos((2(p+q)-n)t)(-1)^{k-q}\int_{\Sigma} p_k d\mu\\
     &+\sum_{0\leq k\leq n,2\nmid k}C_n^k \f{1}{2^n}\sum_{p=0}^{n-k}\sum_{q=0}^{k}C_{n-k}^p C_k^q (-1)^{[\f{k}{2}]}\sin((2(p+q)-n)t)(-1)^{k-q}\int_{\Sigma} p_k d\mu.
\end{align*}
Next let $2(p+q)-n=\pm s$ and sum up in terms of $s$ first. We get
\begin{align}
|\Sigma_t|&=\sum_{s=\f{1-(-1)^n}{2},+2}^n \sum_{p+q=\f{n+s}{2}or\f{n-s}{2}}\sum_{q\leq k\leq n-p,2|k}C_n^k \f{1}{2^n}C_{n-k}^p C_k^q (-1)^{\f{k}{2}+k-q}\cos(st)\int_{\Sigma} p_k d\mu\nonumber\\
     &\quad+\sum_{s=\f{1-(-1)^n}{2},+2}^n \sum_{p+q=\f{n+s}{2}or\f{n-s}{2}}\sum_{q\leq k\leq n-p,2\nmid k}C_n^k \f{1}{2^n}C_{n-k}^p C_k^q (-1)^{[\f{k}{2}]+k-q}\nonumber\\
     &\qquad\times(-1)^{\chi_{\{2(p+q)-n\leq 0\}}}\sin(st)\int_{\Sigma} p_k d\mu\nonumber\\
     &\leq \sum_{s=\f{1-(-1)^n}{2},+2}^n \sqrt{(E(s))^2+(F(s))^2},\label{ineq1}
\end{align}
where
\begin{align*}
E(s)&=\sum_{p+q=\f{n+s}{2}or \f{n-s}{2}}\sum_{q\leq k\leq n-p,2|k}C_n^k \f{1}{2^n}C_{n-k}^p C_k^q (-1)^{\f{k}{2}+k-q}\int_{\Sigma} p_k d\mu,\\
F(s)&=\sum_{p+q=\f{n+s}{2}or \f{n-s}{2}}\sum_{q\leq k\leq n-p,2\nmid k}C_n^k \f{1}{2^n}C_{n-k}^p C_k^q (-1)^{[\f{k}{2}]+k-q}(-1)^{\chi_{\{2(p+q)-n\leq 0\}}}\int_{\Sigma} p_k d\mu,
\end{align*}
and ``$+2$'' means that the step length of the summation for $s$ is $2$.

Next we show that for the geodesic sphere with radius $r\in [0,\f{\pi}{2})$, the equality holds. For this special hypersurface, $\int_{\Sigma} p_k d\mu=\omega_n \sin^n r\cot^k r=\omega_n\sin^{n-k}r \cos^k r$. Thus
\begin{align*}
|\Sigma_t|&=\omega_n\sum_{k=0}^n C_n^k (\cos t \sin r)^{n-k}(\sin t\cos r)^k\\
     &=\omega_n\sin^n(r+t)\\
     &=\omega_n\f{1}{2^n}\sum_{q=0}^nC_n^q \cos ((2q-n)(r+t)-\f{n}{2}\pi)(-1)^{n-q}.
\end{align*}
For simplicity, we assume $n$ is even. Then
\begin{align*}
|\Sigma_t|&=\omega_n\f{1}{2^n}\sum_{q=0}^nC_n^q \cos ((2q-n)(r+t))(-1)^{\f{n}{2}+n-q}\\
     &=\omega_n \sum_{s=0,2,\cdots,n}\sum_{2q-n=\pm s}\f{1}{2^n}C_n^q \cos(s(r+t))(-1)^{\f{3n}{2}-q}\\
     &=\omega_n \sum_{s=0,2,\cdots,n}\sum_{2q-n=s}2\f{1}{2^n}C_n^q \cos(s(r+t))(-1)^{\f{3n}{2}-q},
\end{align*}
where we noted that the coefficients of $\cos(s(r+t))$ for two choices of $q$ are the same.

Now expand $\cos(s(r+t))=\cos sr \cos st-\sin sr \sin st$. We find that all the inequalities in \eqref{ineq1} become equalities for $t=\f{\pi}{2}-r$ and $|\Sigma_t|=\omega_n\sin^n(\f{\pi}{2})=\omega_n$. Thus for the geodesic sphere, the equality in \eqref{ineq1} holds.

On the other hand, assume the equality holds. Then when some $t$ satisfies $|V_t|=V(\f{\pi}{2})=\omega_n\int_0^\f{\pi}{2} \sin^n r dr $, there must hold $|{\Sigma}_t|= S(\f{\pi}{2})=\omega_n$ for this $t$. So the isoperimetric inequality implies that $\Sigma_t=\SS^n(1)$. Then the initial hypersurface must be a geodesic sphere.

Thus Theorem \ref{thm2} is proved.

\begin{rem}\label{rem:5.1}
In section \ref{sec:1}, we discussed the special case $n=2$ of \eqref{eq-thm2}, which is just the Minkowski-type inequality for convex surfaces in $\SS^3$. Here, to get more feeling of the inequality \eqref{eq-thm2}, we give the precise expressions for $n=3$ and $n=4$. For $n=3$, we have
\begin{align*}
\omega_3&\leq \sqrt{(\f{1}{4}(|{\Sigma}|-3\int_{\Sigma} p_2d\mu))^2+(\f{1}{4}(3\int_{\Sigma} p_1d\mu-\int_{\Sigma} p_3d\mu))^2}\\
     &+\sqrt{(\f{3}{4}(|{\Sigma}|+\int_{\Sigma} p_2 d\mu))^2+(\f{3}{4}(\int_{\Sigma} p_1d\mu+\int_{\Sigma} p_3d\mu))^2}.
\end{align*}
And for $n=4$, we have
\begin{align}
\omega_4&\leq \sqrt{(\f{1}{8}(|{\Sigma}|-6\int_{\Sigma} p_2d\mu+\int_{\Sigma} p_4 d\mu))^2+(\f{1}{2}(\int_{\Sigma} p_1 d\mu-\int_{\Sigma} p_3d\mu))^2}\nonumber\\
        &+\sqrt{(\f{1}{2}(|{\Sigma}|-\int_{\Sigma} p_4d\mu))^2+(\int_{\Sigma} p_1d\mu+\int_{\Sigma} p_3d\mu)^2}\nonumber\\
        &+\f{3}{8}(|{\Sigma}|+2\int_{\Sigma} p_2d\mu+\int_{\Sigma} p_4d\mu).\label{eq-sphere-4}
\end{align}
Note that for 4-dimensional hypersurface $\Sigma$ in $\SS^5$, we have the Gauss-Bonnet-Chern formula
\begin{equation}\label{gb-sphere}
  \int_{\Sigma}(p_4+2p_2+1)d\mu=\frac 1{4!}\int_{\Sigma}L_2d\mu=\omega_4.
\end{equation}
Therefore, the inequality \eqref{eq-sphere-4} can be further simplified by using the formula \eqref{gb-sphere}.
\end{rem}

\begin{rem}
As in the hyperbolic case, when the hypersurface $\Sigma\subset\SS^{5}$ is small, the inequality \eqref{eq-sphere-4} reduces to the Euclidean version \eqref{eq-euclid}. This can be seen using the similar argument as in Remark \ref{rem:hyper}.
\end{rem}

\section{Sobolev type inequalities on convex hypersurfaces in sphere}

In this section, by using a different method with the previous sections, we give the proof of Theorem \ref{alex-fenchel-sphere}.

\subsection{Evolution equations}

Considering $\Sigma$ as the initial hypersurface, we flow $\Sigma$ in $\SS^{n+1}$ under the flow equation $X:\Sigma\times[0,T^*)\ra\SS^{n+1}$
\begin{equation*}
 \pt_tX=F\nu,
\end{equation*}
where $F$ is a curvature function and $\nu$ is the unit normal to the flow hypersurfaces $\Sigma_t$. By a standard calculation as in \cite{Hui,GL,LWX}, we have the following evolution equations.

\begin{lem}
Under the curvature flow $\pt_tX=F\nu$ in $\SS^{n+1}$, we have
\begin{eqnarray}
 \frac d{dt}|\Sigma_t|&=&n\int_{\Sigma_t}Fp_1d\mu_t\label{evl-area}\\
\frac d{dt}\int_{\Sigma_t}p_md\mu_t&=&(n-m)\int_{\Sigma_t}Fp_{m+1}d\mu_t-m\int_{\Sigma_t}Fp_{m-1}d\mu_t.\label{evl-pm}
\end{eqnarray}
\end{lem}
\proof
The evolution of \eqref{evl-area} follows by the same calculation as in \cite{Hui}. To obtain \eqref{evl-pm}, we can calculate as in \cite{Hui,GL,LWX}. Firstly, we have as in \cite{Hui}
\begin{eqnarray*}
  \pt_th_i^j &=& -\nabla^j\nabla_iF-Fh_i^kh_k^j-F\delta_i^j.
\end{eqnarray*}
Then
\begin{eqnarray*}
  \pt_t\sigma_m &=& \frac{\pt \sigma_m}{\pt h_i^j}\pt_th_i^j\\
  &=& -\frac{\pt \sigma_m}{\pt h_i^j}\nabla^j\nabla_iF-F\frac{\pt \sigma_m}{\pt h_i^j}h_i^kh_k^j-F\frac{\pt \sigma_m}{\pt h_i^j}\delta_i^j\\
  &=&-\nabla^j((T_{m-1})_j^i\nabla_iF)-F(\sigma_1\sigma_m-(m+1)\sigma_{m+1})-(n+1-m)F\sigma_{m-1},
\end{eqnarray*}
where we used the facts that
\begin{equation*}
  (T_{m-1})_j^i=\frac{\pt \sigma_m}{\pt h_i^j}=\frac 1{(m-1)!}\delta_{j_1\cdots j_{m-1}\,j}^{i_1\cdots i_{m-1}\,i}h_{i_1}^{j_1}\cdots h_{i_{m-1}}^{j_{m-1}},
\end{equation*}
is divergence-free on hypersurface in space forms (see \cite{Rei}), and that
\begin{eqnarray*}
  \frac{\pt \sigma_m}{\pt h_i^j}h_i^kh_k^j&=&\sigma_1\sigma_m-(m+1)\sigma_{m+1}\\
  \frac{\pt \sigma_m}{\pt h_i^j}\delta_i^j&=&(n+1-m)\sigma_{m-1}.
\end{eqnarray*}
Therefore we obtain that
\begin{eqnarray*}
  \frac d{dt}\int_{\Sigma_t}p_md\mu_t &=&  \frac 1{C_n^m}  \frac d{dt}\int_{\Sigma_t}\sigma_md\mu_t\\
  &=&\frac 1{C_n^m}\int_{\Sigma_t}(\pt_t\sigma_md\mu_t+\sigma_m\pt_td\mu_t )\\
  &=&(n-m)\int_{\Sigma_t}Fp_{m+1}d\mu_t-m\int_{\Sigma_t}Fp_{m-1}d\mu_t.
\end{eqnarray*}
\endproof

To simplify the notation, in the following we denote
\begin{eqnarray}\label{tilde-lk}
\tilde{L}_k &=& \frac 1{C_n^{2k}(2k)!}L_k=\,\sum_{i=0}^kC_k^ip_{2k-2i}.
\end{eqnarray}

\begin{lem}\label{evl_L_k}
Under the curvature flow $\pt_tX=F\nu$ in $\SS^{n+1}$, we have
\begin{eqnarray*}
\frac d{dt}\int_{\Sigma_t}\tilde{L}_k d\mu_t&=&(n-2k)\sum_{i=0}^kC_k^i\int_{\Sigma_t}F{p_{2k-2i+1}}d\mu_t.
\end{eqnarray*}
\end{lem}
\proof
The proof is by a direct calculation
\begin{eqnarray*}
\frac d{dt}\int_{\Sigma_t}\tilde{L}_k d\mu_t&=& \sum_{i=0}^kC_k^i\frac d{dt}\int_{\Sigma_t}p_{2k-2i}d\mu_t \\
&=&\sum_{i=0}^kC_k^i\int_{\Sigma_t}\left((n-2k+2i)F{p_{2k-2i+1}}-2(k-i)F{p_{2k-2i-1}}\right)d\mu_t\\
&=&\sum_{i=0}^kC_k^i\int_{\Sigma_t}(n-2k+2i)F{p_{2k-2i+1}}d\mu_t\\
&&\quad -\sum_{i=1}^kC_k^{i-1}\int_{\Sigma_t}2(k-i+1)F{p_{2k-2i+1}}d\mu_t\\
&=&(n-2k)\sum_{i=0}^kC_k^i\int_{\Sigma_t}F{p_{2k-2i+1}}d\mu_t.
\end{eqnarray*}
\endproof

\subsection{Proof of Theorem \ref{alex-fenchel-sphere}}

Recently, Makowski-Scheuer \cite{Mak-Sch} and Gerhardt \cite{gerh} studied the curvature flows in the sphere. If the initial hypersurface $\Sigma\subset\SS^{n+1}$ is closed and strictly convex, then under the inverse mean curvature flow
\begin{equation*}
  \pt_tX=\frac 1{p_1}\nu,
\end{equation*}
there exists a finite time $T^*<\infty$ such that the flow hypersurface $\Sigma_t$ converges to an equator in $\SS^{n+1}$ and the mean curvature of $\Sigma_t$ converges to zero almost everywhere in the sense of (see Theorem 1.4 in \cite{Mak-Sch})
\begin{equation}\label{eq-asym}
  \lim_{t\rightarrow T^*}\int_{\Sigma_t}p_1^{\alpha}d\mu_t=0
\end{equation}
for all $1\leq\alpha<\infty$.

For each $t\in[0,T^*)$, define the quantity $Q(t)$ by
\begin{equation}
  Q(t)=|\Sigma_t|^{-\frac{n-2k}{n}}\int_{\Sigma_t}\tilde{L}_kd\mu_t.
\end{equation}
On the one hand, by Lemma \ref{evl_L_k} and Lemma \ref{newton_ineq} (note that strictly convex implies all principal curvatures of $\Sigma_t$ are positive, and certainly belongs to $\Gamma_k^+$), we have
\begin{eqnarray*}
\frac d{dt}\int_{\Sigma_t}\tilde{L}_k d\mu_t&=&(n-2k)\sum_{i=0}^kC_k^i\int_{\Sigma_t}\frac{p_{2k-2i+1}}{p_1}d\mu_t\\
&\leq &(n-2k)\sum_{i=0}^kC_k^i\int_{\Sigma_t}p_{2k-2i}d\mu_t\\
&=&(n-2k)\int_{\Sigma_t}\tilde{L}_kd\mu_t.
\end{eqnarray*}
Equality holds if and only if $\Sigma_t$ is totally umbilical. On the other hand, the area of the flow hypersurface evolves as
\begin{eqnarray*}
\frac {d}{dt}|\Sigma_t| &=& n|\Sigma_t|.
\end{eqnarray*}
Therefore we obtain that the quantity $Q(t)$ is monotone decreasing in $t$, i.e.,
\begin{equation}\label{monotone-Q}
  \frac d{dt}Q(t)\leq 0.
\end{equation}

Since under the inverse mean curvature flow, the flow hypersurfaces converge to an equator in $\SS^{n+1}$ and the mean curvature of $\Sigma_t$ converges to zero almost everywhere in the sense of \eqref{eq-asym}, we have
\begin{equation}\label{limit-Q}
  \lim_{t\rightarrow T^*}Q(t)=\omega_{n}^{\frac {2k}{n}}.
\end{equation}
Combining \eqref{monotone-Q} and \eqref{limit-Q}, we have
\begin{equation*}
  Q(0)=|\Sigma|^{-\frac{n-2k}n}\int_{\Sigma}\tilde{L}_kd\mu\geq\,\lim_{t\ra T^*}Q(t)=\omega_{n}^{\frac {2k}{n}}.
\end{equation*}
Hence noting \eqref{tilde-lk}, we obtain that
\begin{equation}\label{eq-lk}
  \int_{\Sigma}L_kd\mu\geq C_{n}^{2k}(2k)!\omega_{n}^{\frac {2k}{n}}|\Sigma|^{\frac{n-2k}{n}}.
\end{equation}
Equality holds in \eqref{eq-lk} if and only if $Q(t)$ is constant in $t$. Then $\Sigma_t$ is totally umbilical for each $t\in [0,T^*)$, and in particular $\Sigma$ is totally umbilical and hence a geodesic sphere. The inequality \eqref{eq-lk} says that the induced metric of convex hypersurfaces in $\SS^{n+1}$ satisfy the optimal Sobolev inequalities. See \cite{GWW2} for further information about the Sobolev inequalities of the same type.


\bibliographystyle{Plain}

\begin{thebibliography}{10}

\bibitem{Alex1} A.D. Alexandrov, \emph{Zur Theorie der gemischten Volumina von konvexen K\"{o}rpern, II. Neue Ungleichungen zwischen den gemischten Volumina und ihre Anwendungen}, Mat. Sb. (N.S.) \textbf{2} (1937) 1205--1238 (in Russian).
\bibitem{Alex2} A.D. Alexandrov, \emph{Zur Theorie der gemischten Volumina von konvexen K\"{o}rpern, III. Die Erweiterung zweeier Lehrsatze Minkowskis \"{u}ber die konvexen Polyeder auf beliebige konvexe Flachen}, Mat. Sb. (N.S.) \textbf{3} (1938) 27--46 (in Russian).

\bibitem{ben} B. Andrews, \emph{Aleksandrov-Fenchel inequalities and curvature flows}, Centre for Mathematics and its
Applications, School of Mathematical Sciences, 1993.

\bibitem{BHW} S. Brendle, P.-K. Hung, and M.-T. Wang, {\it A Minkowski-type inequality for hypersurfaces in the Anti-deSitter-Schwarzschild manifold}, arXiv: 1209.0669.

\bibitem{Bla} W. Blaschke, \emph{\"{U}ber eine geometrische Frage von Euklid bis heute}, Hamburger Mathematiches
Einzelschriften, vol. \textbf{23}, 1938.

\bibitem{BM} A. A. Borisenko and V. Miquel, {\it Total curvatures of convex hypersurfaces in hyperbolic space}, Illinois J. Math. \textbf{43} (1999), 61--78.

\bibitem{Chern}S. S. Chern, \emph{A simple intrinsic proof of the Gauss-Bonnet formula for closed Riemannian manifolds}. The Annals of Mathematics, 1944, \textbf{45}(4): 747--752.

\bibitem{CW} S.-Y. A. Chang and Y.Wang, {\it On Aleksandrov-Fenchel inequalities for k-convex domains}, Milan J. Math., \textbf{79} (2011), no. 1, 13--38.

\bibitem{de-girao} L.L. de lima and F. Gir\~{a}o, {\it An Alexandrov-Fenchel-type inequality in hyperbolic space with an application to a Penrose inequality}, arXiv:1209.0438.

\bibitem{do-carmo}M. do Carmo and F. Warner, {\it Rigidity and convexity of hypersurfaces in spheres}, J. Differ.
Geom. \textbf{4} (1970), 133--144.

\bibitem{Fe}W. Fenchel, \emph{In\'{e}galit\'{e}s quadratiques entre les volumes mixtes des corps convexes}, C. R. Acad. Sci. Paris S\'{e}r. I Math. \textbf{203} (1936), 647--650.

\bibitem{gerh}C. Gerhardt, {\it Curvature flows in the sphere}, arXiv:1308.1607


\bibitem{guan}P. Guan, {\it Topics in Geometric Fully Nonlinear Equations}, available at \url{http://www.math.mcgill.ca/guan/notes.html}.


\bibitem{GL} P. Guan and J. Li, {\it The quermassintegral inequalities for k-convex starshaped domains}, Adv. Math. \textbf{221 }(2009), 1725--1732.

\bibitem{GS} E. Gallego and G. Solanes, {\it Integral geometry and geometric inequalities in hyperbolic space}, Differential Geom. Appl. 22 (2005), 315--325.

\bibitem{GWW0} Y. Ge, G. Wang and J. Wu, \emph{A new mass for asymptotically flat manifolds}, arXiv:1211.3645.

\bibitem{GWW1} Y. Ge, G. Wang and J. Wu, \emph{Hyperbolic Alexandrov-Fenchel quermassintegral inequalities I}, arXiv:1303.1714.

\bibitem{GWW2}Y. Ge, G. Wang and J. Wu, \emph{Hyperbolic Alexandrov-Fenchel quermassintegral inequalities II}, arXiv:1304.1417.

\bibitem{GWW3}Y. Ge, G. Wang and J. Wu, \emph{The GBC mass for asymptotically hyperbolic manifolds}, arXiv:1306.4233.

\bibitem{Hui}G. Huisken, \emph{Flow by mean curvature of convex surfaces into spheres}, J. Differential Geom., \textbf{20}(1984),237--266.

\bibitem{Kno} H. Knothe, \emph{Zur Theorie der konvexen K\"{o}rper im raum konstanter positiver Kr\"{u}mmung,
Revista da Faculdade de Ci\^{e}ncias de Lisboa}, S\'{e}rie 2 A \textbf{2} (1952), 336--348.

\bibitem{LWX} H. Li, Y.Wei and C. Xiong, \emph{A geometric inequality on hypersurface in hyperbolic space}, arXiv:1211.4109.

\bibitem{Mink} H. Minkowski, \emph{Volumen und Oberfl\"{a}che}, Math. Ann. \textbf{57} (1903), 447--495.

\bibitem{McC} J. McCoy, \emph{Mixed volume preserving curvature flows}, Calc. Var. \textbf{24}, (2005), 131--154.

\bibitem{MR}  S. Montiel and A. Ros,  {\it Compact hypersurfaces: the Alexandrov theorem for higher order mean curvatures}, Differential Geometry, Blaine Lawson and Keti Tonenblat, Pitman Monographs, 1991, \textbf{52}: 279--297.


\bibitem{Mak-Sch} M. Makowski and J. Scheuer, {\it Rigidity results, inverse curvature flows and Alexandrov-Fenchel type inequalities in the sphere}, arXiv:1307.5764.

\bibitem{Na} J. Nat\'{a}rio, {\it A Minkowski-type inequality for convex surfaces in the hyperbolic 3-space}, arXiv:1307.4239.

\bibitem{Oss} R. Osserman, \emph{The isoperimetric inequality}, Bull. Amer. Math. Soc. \textbf{84} (1978), 1182--1238.

\bibitem{Ros} A. Ros, \emph{The isoperimetric problem}, Lecture series at the Clay Mathematics Institute Summer School on the Global Theory of Minimal Surfaces, summer 2001,  Mathematical Sciences Research Institute, Berkeley, California, available at \url{http://www.ugr.es/~aros/isoper.pdf}

\bibitem{Rei}R. Reilly, \emph{On the Hessian of a function and the curvatures of its graph}, Michigan Math. J., \textbf{20}(1973),373--383.

\bibitem{S} E. Schmidt, {\it Die isoperimetrischen Ungleichungen auf der gew\"{o}hnlichen Kugel und f\"{u}r Rotationsk\"{o}rperim n-dimensionalen sph\"{a}rischen Raum}, (German) Math. Z. \textbf{46}, (1940), 743--794.

\bibitem{Sant}L. Santal\'{o}, \emph{A relation between mean curvatures of parallel convex bodies in spaces of constant curvature}, Rev. Un. Mat. Argentina \textbf{21} (1963), 131--137.

\bibitem{Schn}R. Schneider, \emph{Convex bodies: the Brunn-Minkowski theory}, Cambridge University Press, 1993.

\bibitem{WX} G. Wang and C. Xia, {\it Isoperimetric type problems and Alexandrov-Fenchel type inequalities in the hyperbolic space}, arXiv:1304.1674.


\end{thebibliography}

\end{document}